\documentclass[11pt]{article}

\usepackage{pdfsync}
\usepackage{amsmath}
\usepackage{amsfonts}
\usepackage{amsthm}
\usepackage{nicefrac}
\usepackage{dsfont}
\usepackage{amssymb}
\usepackage{mathrsfs}
 \usepackage{graphicx}
\usepackage{xcolor, caption, hyperref, psfrag}

\newtheorem{lm}{Lemma}
\newtheorem{te}{Theorem}
\newtheorem{cor}{Corollary}
\newtheorem{re}{Remark}
\newtheorem{prop}{Proposition}
\newtheorem{con}{Conjecture}

\def\bl{\begin{lm}}
\def\el{\end{lm}}
\def\bpr{\begin{prop}}
\def\epr{\end{prop}}
\def\bt{\begin{te}}
\def\et{\end{te}}
\def\bc{\begin{cor}}
\def\ec{\end{cor}}
\def\bcon{\begin{con}}
\def\econ{\end{con}}
\def\br{\begin{re}}
\def\er{\end{re}}
\def\be{\begin{equation}}
\def\ee{\end{equation}}
\def\bp{\begin{proof}}
\def\ep{\end{proof}}

\def\RR{{\mathbb R}}

\def\PP{{\mathbb P}}
\def\DD{{\mathbb D}}

\def\NN{{\mathbb N}}

\def\MV{{\mathcal V}}
\def\MG{{\mathcal G}}

\def\MR{{\mathcal R}}

\def\g{{\gamma}}
\def\t{{\tau}}

\def\1{{\mathbf 1}}

\allowdisplaybreaks

\begin{document}
\sloppy

\title {The Kingman tree length process has infinite quadratic variation\footnote{Work partially supported by the DFG Priority Programme SPP 1590 ``Probabilistic Structures in Evolution''.}}

\author{ Iulia Dahmer
\thanks{Institut f\"ur Mathematik, Goethe-Universit\"at, 60054 Frankfurt am Main, Germany. \newline\texttt{dahmer@math.uni-frankfurt.de, wakolbinger@math.uni-frankfurt.de}
} \qquad
 Robert Knobloch\thanks{Fachrichtung Mathematik, Universit\"at des Saarlandes,  Postfach 151150, 66041 Saarbr\"ucken, Germany.
\newline
\texttt{knobloch@math.uni-sb.de}}
 \qquad
Anton Wakolbinger$^\dagger$%
}
\date{\today}
\maketitle 
\begin{abstract}
In the case of neutral populations of fixed sizes in equilibrium whose genealogies are described by the Kingman $N$-coalescent back from time $t$ consider the associated processes of total tree length as $t$ increases. We show that the (c\`adl\`ag) process to which the sequence of compensated tree length processes converges as  $N$ tends to infinity is a process of infinite quadratic variation; therefore this process cannot be a semimartingale. This answers a question posed in \mbox{Pfaffelhuber et al. (2011)}.
\end{abstract}
\begin{small}
\end{small}

\section{Introduction and main result}

The Kingman coalescent is a classical model in mathematical population genetics used for describing the genealogies for a wide class of population models (see e.g \cite{Wa08}). The population models in question are neutral, exchangeable and with an offspring distribution of finite variation. One particular example is the Moran model (\cite{Mo58}). This is a stationary continuous-time model for populations of fixed size $N$ in which the reproduction takes place according to the following rule: starting with a population of size $N$, after an exponential time of parameter $\binom N 2$ a pair of individuals is picked uniformly at random from the population, out of which one individual dies and the other one gives birth to one child.

The ancestry of a Moran population of size $N$ started at time $-\infty$ is at any time $t \in \RR$ described by the Kingman $N$-coalescent. This is a process with values in the set of partitions of $\{1, \dots, N\}$ which starts in the partition in singletons and has the following dynamics (backwards in time): given the  process is in state $\pi_k$, it jumps at rate $\binom k 2$ to a state $\pi_{k-1}$ which is obtained by merging two randomly chosen elements of $\pi_k$. The process can be represented graphically as a binary rooted tree which, when traced back from its $N$ leaves (and correspondingly $N$ {\em external branches}), exhibits a binary merger at rate $\binom k 2$ while there are $k$ branches left.

One particular feature of coalescent trees that has been intensively investigated in the literature, due to its relevance in statistical studies of genetic data, is their total length (the sum of the lengths of all the branches of the tree). In the case of the Kingman coalescent tree started with $N$ leaves the total length is in expectation equal to twice the harmonic number $h_{N-1}=\sum_{i=1}^{N-1}\frac 1 i$ and when $N$ tends to infinity (half of) the total length compensated by  $\log N$ converges in law to a Gumbel distributed random variable. In the case of coalescent processes with multiple mergers the total length has been studied in various papers, for instance \cite{bebeli}, 
\cite{bebesw}, \cite{bebesw2}, \cite{dr}, \cite{ikmo}, \cite{ke}, \cite{mo}. 

As time $t$ increases the Moran population evolves and its genealogy changes, giving rise to a tree-valued process $\mathcal R^N=(\mathcal R^N_t)_{t \in \RR}$, the {\it evolving Kingman $N$-coalescent}. The associated process of total tree length was investigated in \cite{PWW11}. (See also the more recent papers \cite{KSW14} and \cite{S12} on the evolution of the total length in the multiple merger case.) Let $\ell(\mathcal R^N_t)$ denote the length of the tree $\mathcal R^N_t$ and call
\[
 \mathfrak L^N=\Big(\mathfrak L^N_t\Big)_{t \in \RR}:=\left(\ell(\mathcal R^N_t)-2\log N\right)_{t \in \RR} 
\]
the {\it compensated tree length process}. Pfaffelhuber et al. \cite{PWW11} investigated the asymptotic behaviour of this process as the population size $N \to \infty$ and showed that there exists a process $\mathfrak L=(\mathfrak L_t)_{t \in \RR}$ with sample paths in $\DD$, the space of c\`adl\`ag functions equipped with the Skorokhod topology, such that 
\begin{equation}\label{weak}
\mathfrak L^N \to \mathfrak L \text{ in law as } N \to \infty.
\end{equation}
The process $\mathfrak L$ is the {\em Kingman tree length process}.
\begin{center}
\includegraphics[scale=0.65]{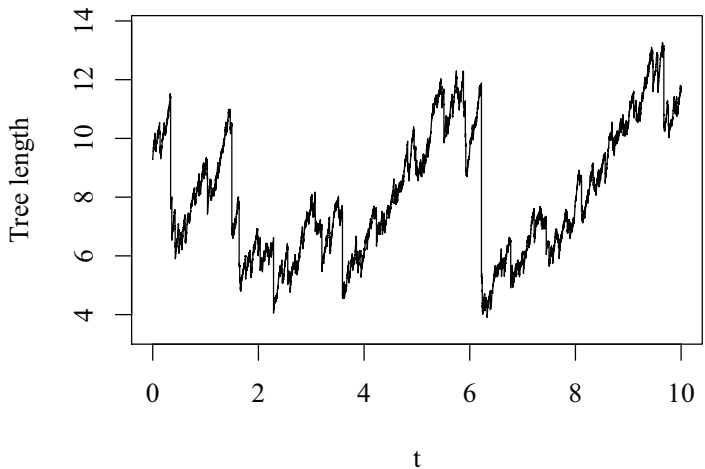}
\captionof {figure}{A realization of the compensated tree length process $\mathfrak L^{N}$ for $N=30$ (courtesy of Peter Pfaffelhuber)}
\label{Fig_tree}
\end{center}
The weak convergence \eqref{weak} can be lifted to convergence in probability, provided a representation for Moran populations on the same probability space for all population sizes $N \in \NN$ is considered. Such a representation is given by the {\it look-down} construction of Donnelly and Kurtz (\cite{DK96}, \cite{DK99}) which  encodes the evolving coalescent in a path-wise consistent way for increasing $N$. If $\mathfrak L^{ld,N}_t$ denotes the compensated length of the tree at time $t$ in the look-down representation, and $\mathfrak L^{ld,N} := (\mathfrak L^{ld,N}_t)_{t\in \mathbb R}$, then, as shown in \cite{PWW11} Proposition 3.2, there exists a process $\mathfrak L^{ld}$, having the same distribution as $\mathfrak L$, such that 
$$d_{Sk}(\mathfrak L^{ld,N}, \mathfrak L^{ld})\longrightarrow 0$$
holds in probability as  $N \to \infty$, where $d_{Sk}$ denotes the Skorokhod metric. The proof of Proposition~3.2 in \cite{PWW11}  is based on the equality in law of the processes $\mathfrak L^{ld,N}$ and $\mathfrak L^{N}$. In Sec. 2 below we include an argument why this equality in law is valid.

The question that we address in this paper is one formulated in \cite{PWW11}, namely whether $\mathfrak L$ is a semimartingale (i.e. whether it can be written as a sum of a local martingale and a process of locally finite variation that are both adapted to the same filtration), and thus would be an instance  for the classical tools of stochastic analysis.  A necessary condition for a c\`adl\`ag process to be a semimartingale is that its quadratic variation is a.s. finite, see e.g.  \cite{Pr04} Theorem~II.22. In \cite{PWW11}  it was proved that the process $\mathfrak L$ has ``infinite infinitesimal variance'', more precisely, $\tfrac 1{\varepsilon |\log \varepsilon|}\mathbb E[(\mathfrak L_\varepsilon-\mathfrak L_0)^2]\to 4$ as $\varepsilon\to 0$.  This implies that the squared increments $(\mathcal L_{t+\varepsilon}-\mathcal L_t)^2$ are for small $\varepsilon$ (at least in expectation) of a larger order than $\varepsilon$, which suggests that $\mathcal L$ should not have finite quadratic variation. We will 
show that indeed $\mathfrak L$ has a.s. infinite quadratic variation (and hence cannot be a semimartingale).   This will be achieved by investigating the jumps of the process~$\mathfrak L^{ld}$. 


Let us now give a brief description of the look-down construction and explain heuristically our approach. A formal description of the look-down graph will then be given in the next section. The main idea behind the look-down representation is to label the individuals in the population according to the persistence (or longevity) of their offspring: label 1 for the individual with the most persistent progeny, 2 for the second and so on. 

\begin{center}
\psfrag{8}{$\mathfrak L^{ld,5}$}
\includegraphics[scale=0.25]{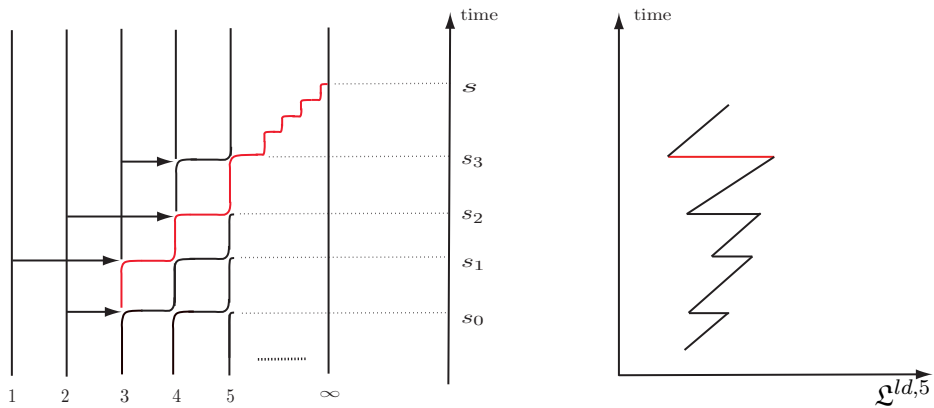}
\captionof {figure}{Detail of one realization of the infinite look-down graph. The line $G$ marked in red is born at time $s_0$ at level 3 and is pushed up one level whenever an arrow is shot towards a level lower than the current level of the line (at times $s_1, s_2, \dots$). In the $N$-look-down graph with $N=5$ the line dies at time $s_3$, whereas in the infinite look-down graph it dies at time $s$ when it reaches level $\infty$. The life-length of the line is equal to $s-s_0$ and its life-length up to $N$ is $s_3-s_0$. The tree length process $\mathfrak L^{ld,5}$ has jumps at the times lines exit level 5. The sizes of these jumps are equal to the lengths of the corresponding lines.}
\label{Fig_ld}
\end{center}

We consider a system of countably many particles describing the sample genealogies ordered by persistence. At any time, each {\em level} 1,2,\ldots is occupied by precisely one particle, and the system evolves as follows: for every pair $i<k$ at rate 1 the particle currently at level $i$ shoots an arrow towards level $k$, independently of everything else. At this time the particle at level $i$ gives birth to a new particle which is placed at level $k$, while for each $j \geq k$ the particle located at level $j$ changes its level from $j$ to $j+1$. To each birth event  we associate a {\it line} which records the time evolution of the levels occupied by the new-born particle (see the graphical representation in Figure \ref{Fig_ld}). This line is pushed up to the next level each time a birth event happens on a level to the left of the current level of the line. We say that the line ends (dies) at the time it reaches level $\infty$. The countable system of all the lines (including the immortal line that sits at 
level 1) makes up the look-down graph with infinitely many levels (or {\em infinite look-down graph} for short); the corresponding representation for a particle system of finite size $N$ is obtained by projecting the infinite look-down graph onto the first $N$ levels. When considering the system with $N$ particles only,  we say that a line dies when it is pushed out of level $N$. Like in the case of the Moran model, the realizations of the sample genealogy can be read off from the look-down graph.

For a line $G$ in the infinite look-down graph we denote by $T^{G}$ its {\it life-length}, i.e. the time span between the birth and the death time of $G$. If we restrict the graph to its first $N$ levels, then $T^{G,N}$, the {\it life-length of the line up to $N$} will denote the time span the line needs until it exits level $N$. In terms of trees, the life-length of a line that dies at some time $t$ in the look-down graph with $N$ levels corresponds to the length of the external branch that falls off the genealogical tree at time $t$. Therefore, the jumps of the (compensated) tree length process $\mathfrak L^{ld,N}$ happen at the times lines exist level $N$ in the $N$-look-down graph and they have sizes equal to the life-lengths up to $N$ of these lines. Hence for $s<t$ we can write
\[\mathfrak L^{ld,N}_t-\mathfrak L^{ld,N}_s =N(t-s) - \sum_G T^{G,N},\]
where the sum is taken over all lines $G$ that exit level $N$ in the time interval $(s,t]$. It was proved in \cite{PWW11} (see Proposition 3.1 and the proof of Proposition 6.1 therein) that for any fixed times $s<t$
\begin{equation}\label{jump_N}
\mathfrak L^{ld,N}_t-\mathfrak L^{ld,N}_s  \longrightarrow \mathfrak L^{ld}_t-\mathfrak L^{ld}_s
\end{equation}
holds in $L^2$, and therefore almost surely along a subsequence $(N_k)_{k \in \NN}$. 

Let us now consider the lines in the infinite look-down graph that die in the time interval $(s,t]$. For every such line there exists an $N$ such that for all $N' \ge N$ this line exits from level $N'$ in the time interval $(s,t]$. Conversely, for any line that exits at level $\infty$ in the complement of the time interval $(s,t]$ there exists an $N$ such for all $N' \ge N$ that this line does not exit from level $N'$ in the time interval $(s,t]$. Therefore, with probability one, it is the life-lengths up to $N_k$ of precisely those lines that reach level $\infty$ in $(s,t]$, which appear as summands on the right-hand side of (\ref{jump_N}) for large enough $k$, and thus contribute to the limit $\mathfrak L_t^{ld, N_k}-\mathfrak L_s^{ld, N_k}$ as $N_k \to \infty$. 

Therefore, in order to understand the jumps of the limiting process $\mathfrak L^{ld}$ that occur in $(s,t]$ one key issue is to understand the behaviour of the life-lengths of the lines that die in the infinite look-down graph in this time interval. The following theorem on the squared life-lengths of these lines is the central ingredient for proving our main result, which is stated in Theorem \ref{main_th} below.

\bt \label{sec_th}
For any $s<t$ the sum of the squared life-lengths of the lines that die in the time interval $(s,t]$ in the infinite look-down graph is almost surely infinite.
\et
\bt \label{main_th}
The Kingman tree length process $\mathfrak L$ has a.s. infinite quadratic variation. That is to say, for any $s<t$ and each sequence $(\mathcal P_n)_{n \in \NN}=\Big((\rho_j^{(n)})_{j=0, \dots,l^{(n)}}\Big)_{n \in \NN}$ of partitions of $[s,t]$ with mesh size tending to zero as $n \to \infty $ one has $\displaystyle \lim_{n \to\infty}\sum_{j=1}^{l^{(n)}} \Big(\mathfrak L_{\rho_j^{(n)}}-\mathfrak L_{\rho_{j-1}^{(n)}}\Big)^2 =\infty$ a.s.
\et
\noindent
We will prove Theorem \ref{main_th} for $\mathfrak L^{ld}$ in place of $\mathfrak L$. This is sufficient, since $\mathfrak L^{ld}$ and $\mathfrak L$ are equal in law.

A key ingredient in the proof of Theorem \ref{sec_th} is the proposition stated below. This result is also of interest in its own right since it sheds light on the overall structure of the look-down graph and the large amount of independence which is built into it. From the brief description of the look-down graph given above (and from the formal definition provided in the next section) it is immediate that  the birth times of lines on some level $k\geq 2$ in the look-down graph form a Poisson process with rate $k-1$. It turns out that the death times of these lines are also points of a Poisson process with the same rate. For the particular case $k=2$ two different proofs of this result were given in \cite{DK06}  and \cite{PW06}.
\bpr\label{prop_poisson}
In the infinite look-down graph, for every $k \in \NN$, $k \geq 2$ consider the process $\eta_k$ of time points at which the lines that were born at level $k$  reach level $\infty$. The processes $\eta_k$ are mutually independent Poisson with rate $k-1$. 
\epr
For each $k=2,3,\ldots$ the process $\mathfrak L^{ld}$ has a jump in each of the points of $\eta_k$.  The size $h$ of this jump is equal to the life-length $T^G$ of the line $G$ that dies at this time point (see the proof of Theorem \ref{main_th}). Let us  emphasize that even though the jump times of $\mathfrak L^{ld}$ are independent,  $\mathfrak L^{ld}$ is not a L\'evy process, because there are dependencies in the jump sizes. Moreover, the integrability condition  $\int_{[0,1]} h^2\nu(dh)<\infty$, which must be satisfied by a L\'evy measure, is violated by the jump intensity measure of $\mathfrak L^{ld}$. Indeed, the expectation of the life-length $T^G$ of a line born at level $k$ is $2/k$ (see \eqref{def_T_G} below) and for large $k$ the distribution of $T^G$ is concentrated around $2/k$ (see the proof of Theorem \ref{sec_th}, which uses a result of \cite{DPS13}). Since the points of $\eta_k$ come at rate $k-1$, the jump intensity measure of $\mathfrak L^{ld}$ has (for large $k$) mass $k-1$ concentrated
around $2/k$.  As a matter of fact, part of the strategy of the proof
of our main result reflects in the simple fact that $\sum (k-1) (\frac 2 k)^2 = \infty$.

\section{The look-down process}

The {\it look-down} construction of Donnelly and Kurtz (\cite{DK96}, \cite{DK99}) is an alternative way of representing the evolution of  Moran (and more general exchangeable) populations, which proves to be a very powerful instrument in investigating population dynamics. As already mentioned in the introduction, this representation of populations of sizes $N$ is done on one and the same probability space for all $N \in \NN$ in such a way that the path-wise consistency of the genealogies is ensured as $N \to \infty$. 

The main idea of the look-down representation is the labeling of the individuals according to the persistence of their offspring in the population. In the first paper \cite{DK96} the persistence of the offspring is taken to hold in probability, whereas in the "modified" look-down construction introduced in \cite{DK99}, this holds almost surely. We will use this second version of the model which we describe below following \cite{PW06}.

We consider a population of infinite size and denote by $\MV$ the set $\RR \times \NN$. An element $(s,i)$ in $\MV$ denotes the individual that occupies level $i$ at time $s$. The levels represent indices given to the individuals in the population according to the persistence of their offspring in the following way: the offspring of the individual that lives at time $s$ at level $i$ almost surely outlives the offspring of any other individual alive at time $s$  on a level $k>i$. The process evolves as follows: to every pair of levels $i, k \in \NN$ with $i<k$ we attach a (rate one) Poisson point process on $\RR$ which we denote by $C_{ik}$. All these Poisson point processes are independent. Each time the clock $C_{ik}$ rings, level $k$  {\em looks down} to level $i$, that is, the current individual at level $i$ reproduces and its offspring is placed at level $k$. For $k \geq 2$ and $s_0 \in \bigcup_{i <k} C_{ik}$ we associate with the individual born at time $s_0$ at level $k$ the set of points
$$G=\bigcup_{j \in \NN_0}[s_j,s_{j+1}) \times \{k+j\},$$
where $s_j:=\inf{\left\{ s>s_{j-1}: s \in \bigcup_{ l<m< k+j}C_{lm}\right\}}$ for $j\in\mathbb N$.
We call $G$ the {\it line} born at time $s_0$ at level $k$ and say that at time $s_j$ the line is pushed from level $k+j-1$ to level $k+j$. Note that a line is pushed one level upwards every time one of the Poisson point processes associated with levels smaller than or equal to the current level of the line experiences an event.  Lines are born on a level $k>1$ at the times of a Poisson point process with rate $(k-1)$ and a line at level $k$ is pushed up with rate $\binom k 2$ because there are $\binom k 2$ independent (rate one) Poisson point processes which trigger the look-down between the levels that are smaller than or equal to $k$. 

We say that a line dies when it reaches level infinity and denote the death time of line $G$ by 
\[d^G:=\lim_{j \to \infty} s_j.\]
Since the rate at which a line is born at a level bigger than or equal to 2 is pushed up is quadratical, we conclude that the time it takes for a line to die is finite almost surely. Level $1$ is never hit by arrows and therefore the offspring of the individuals living on this level persist forever in the population. We call the line $\RR \times \{1\}$ {\it the immortal line}.
\par
The set of all the lines is countable and it forms a partition of $\MV$. The random graph obtained in this way is called the {\it look-down graph} (with infinitely many levels).
This graph records the evolution of a population of infinite size.  Embedded in the look-down process are all the $N$-particle look-down processes corresponding to populations of sizes $N \in \NN$. The $N$-particle look-down process is constructed in a similar way, but the graph has only $N$ levels and we say that a line dies when it exits level $N$. Any $N$-particle look-down process can be recovered as the projection of the infinite look-down process on the first $N$ levels.

The ordering by persistence (corresponding to the direction of the arrows from left to right in Figure \ref{Fig_ld}) induces an asymmetry  in the look-down graph: the offspring size of an individual with a lower level tends to be larger than that of an individual with a larger level. Nevertheless,
the ancestral process back from a fixed time $t$ that is induced by the random look-down graph is the Kingman coalescent. In order to see this,  consider two lines $G$ and $G'$. For $(s,l) \in G$ and $(t,i) \in G'$ with $s \le t$ we say that $(s,l)$ is the {\it ancestor} of $(t,i)$ and we write
\[A_s(t,i)=l,\]
if either the two lines are the same or there are some lines $G_1,\dots,G_m$ such that $G_1$ descends from $G$, $G_k$ descends from $G_{k-1}$, for $k=2,\dots,m$ and $G'$ descends from $G_m$. Two individuals $(t,i)$ and $(t,j)$ living at time $t$ have the same ancestor at time $s$ if $A_s(t,i)=A_s(t,j)$ and we write $i\stackrel {u}{\sim} j$ with $u=t-s$. The random equivalence relation $\stackrel {u}{\sim}$ defines the ancestral process of the population alive at time $t$, $ \MR^{ld}_t=( \MR^{ld}_t(u))_{u\in \RR}$. It is not difficult to check that for each $t$ the restriction $\MR^{ld,N}_t$ of $\MR^{ld}_t$ to $\{1, \dots, N\}$  is equal in law to the $N$-Kingman coalescent $\mathcal R_t^N$, when both are viewed as metric trees. The consistency property then implies that the genealogy $\MR^{ld}_t$ of the infinite population has the distribution of the Kingman coalescent  $\MR_t$.

%

The trees $\MR^{ld,N}_t$ and $\MR^N_t$ come with a labeling of their leaves by $1,\ldots, N$, which in the case of $\MR^{ld,N}_t$ corresponds to the levels. It is important to note that, for $N \in \mathbb N$, the tree length processes $\mathfrak L^{ld,N} $ and $ \mathfrak L^N $ have the same distribution, even though for $N\ge 3$  the distributions of the leaf-labeled metric tree-valued processes $\MR^{ld,N}$ and $\MR^N$  are different.  As already stated above, for any fixed time $t$, the distribution of  $\MR^{ld,N}_t$ equals that of  $\MR^{N}_t$. Moreover, this distribution is {\em exchangeable}, i.e. invariant under a permutation of the labels.
If one considers instead of the leaf-labeled trees the {\em unlabeled trees} (i.e. the equivalence classes of leaf-labeled trees  under all permutations of the labeling), then it is clear how the Moran dynamics acts on these unlabeled trees: after an  exponential time with parameter $N\choose 2$ (at time $\tau$, say) a pair of leaves is chosen completely at random, one of them to die, which results in the removal of the external branch that is below the leaf that dies at time~$\tau$, and the other
to be parental, which results in two leaves having distance $0$ at time $\tau$.  (A~formal description of this so-called {\em tree-valued Moran dynamics of population size $N$} is given in \cite{GPW13}, Def.~2.18.) With the look-down dynamics acting on the trees whose leaves are labeled by the levels, it is always the leaf at level $N$ that dies, and compared to the Moran dynamics there is a bias towards the lower levels in the choice of the parental leaf. However, because the distribution of  $\MR^{ld,N}_{\tau-}$ (like that of $\MR^{ld,N}_t$) is invariant under permutations of the labels, the choice from the labeled leaves of $\MR^{ld,N}_{\tau-}$ (in spite of its bias) amounts to a uniform choice of a pair of leaves from the unlabeled tree that corresponds to $\MR^{ld,N}_{\tau-}$.  (More formally, for a leaf-labeled tree $x$, denote the unlabeled tree obtained from $x$ by $\Phi(x)$, and write $P(x,.)$ for the look-down transition probability in one reproduction step starting from $x$. Also, for an unlabeled 
tree $y$, denote by $\Lambda(y,.)$ the uniform distribution on the $n!$ leaf-labeled trees in the equivalence class described by $y$, and write $Q(y,.)$ for the Moran transition probability on the unlabeled trees  in one reproduction step starting from $x$. What we have just explained amounts to the relation $\Lambda  P = Q \Lambda$, which is one of the two criteria in Theorem 2 of \cite{RP81}. The other criterion in this theorem (requiring that $\Lambda \Phi = I$) is clearly satisfied. Hence, this theorem yields (first for the chains embedded at the reproduction times and then also for the processes in   continuous time) that $\Phi(\MR^{ld,N})$ is a Markov process whose transitions are given by $Q$.) Altogether, this shows that the lookdown dynamics yields the same Markovian projection on the unlabeled trees as the Moran dynamics. Since the tree length is a functional of the unlabeled tree, this shows that $\mathfrak L^{ld,N} $ and $ \mathfrak L^N $ have the same distribution.

With a view towards the jumps of $\mathfrak L^{ld}$, in the look-down graph with infinitely many levels let us consider a line $G$ born at level $l^G \geq 2$. The time this line needs in order to reach level infinity is
\be\label{def_T_G}
T^G=\sum_{j=l^G}^{\infty} X_j^G, 
\ee
where the time $X_j^G$ spent by the line at level $j$ is an exponentially distributed random variable with parameter $\binom j 2$ and the $X_j$'s are independent from one another for different $j$'s. We call $T^G$ the {\it life-length} of the line $G$. In terms of trees, the life-length of a line that dies at some time $t$ represents the length of the external branch that falls off the genealogical tree at time $t$. When restricting to the first $N$ levels in the graph, we define 
\be\label{def_T_GN}
T^{G,N}:=\sum_{j=l^G}^N X_j^G,
\ee
to be the {\it  life-length up to level $N$} of the line $G$.
\section{Proof of Theorem~\ref{sec_th}}
Before we embark on proving Theorem~\ref{sec_th} let us
provide the proof of Proposition \ref{prop_poisson} which is a key ingredient in the
proof of this theorem.

{\it Proof of Proposition \ref{prop_poisson}}

For every $n \ge 2$ and $2\leq k \leq n$ let us write $\eta_k^n$ for the process of arrival times at level $n$ of lines born at level $k$. For $k=n$, the process $\eta_n^n$ equals the process of time points were new lines are born at level $n$. Since for each $1\leq m \leq n-1$ new lines at level $n$ are born  via birth events triggered from level $m$ at rate $1$, independently of everything else, it is clear that for every $n \ge 2$, $\eta_n^n$ is a Poisson process with rate $n-1$ that is independent of $(\eta^n_2, \ldots, \eta^n_{n-1})$.

It is thus sufficient to prove the following claim: \\\\
\indent $(\ast)$ \,  for every $n \geq 2$ the processes $\eta_k^n$, $2\leq k \leq n-1$,  are Poisson processes of rate $k-1$ and they are independent from one another for $n$ fixed and different values of $k$.

Assuming this claim  holds, remember that for a line 
\[G=\bigcup_{j \in \NN_0}[s^G_j,s^G_{j+1}) \times \{k+j\}\]
born at level $k$, the time point $s^G_j$ is the time the line reaches level $k+j$ and that its death time
\[
d^G=\lim_{j \to \infty} s^G_j
\]
is finite almost surely. Now, denoting by $\MG_k$ the set of all the lines in the look-down graph which are born at level $k$, it follows that the time points $\{s^G_j\}_{G \in \MG_k}$ are the  points of the process $\eta_k^{k+j}$, whereas the points $\{d^G\}_{G \in \MG_k}$ are the  points of the process $\eta_k$. Thus, the assertion of the proposition follows from the claim.

We now prove the claim $(\ast)$ by an induction argument.

For the basic step of the induction let $n=3$. At level 2, lines are born at the times of the Poisson process $\eta_2^2$ and every time a line is born, the line that occupied the level 2 is pushed up to level 3. Therefore, a line born at level 2 arrives at level 3 at the next time point of $\eta_2^2$ after the line's birth time. It follows that the set of points of $\eta_2^2$ is equal to the set of points of $\eta_2^3$ and hence $\eta_2^3$ is a Poisson process with rate 1. Moreover, $\eta_2^3$ and $\eta_3^3$  are independent.
\begin{flushleft}
\begin{minipage}{12.5cm}
$\quad$ We assume now that the claim holds for $n$ and prove it for $n+1$. From the induction assumption and the last sentence in the first paragraph of this proof it follows that the processes $\eta_k^n$, $k=2,\ldots, n,$  are independent Poisson processes of rate $k-1$. A fortiori,  the process $\eta^n_{2, \dots, n}$ of arrival times at level $n$, obtained by superposing the independent processes $\eta_2^n, \dots, \eta_n^n$ is Poisson with rate $\binom n 2$.  A line currently at level $n$ is pushed to level $n+1$ at the next point of $\eta^n_{2, \dots, n}$ after the line's arrival at level $n$. Therefore, there is a bijective function $\phi$ from the collection of points of the process $\eta^n_{2, \dots, n}$ into itself which maps the time a line arrives at level $n$ onto the time it is pushed up (and arrives at level $n+1$) (see Figure \ref{Fig_eta}). 
\\
\end{minipage}
\begin{minipage}{3.5cm}
\begin{center}
\psfrag{a}{$\phi$}
\includegraphics[height=5.5cm]{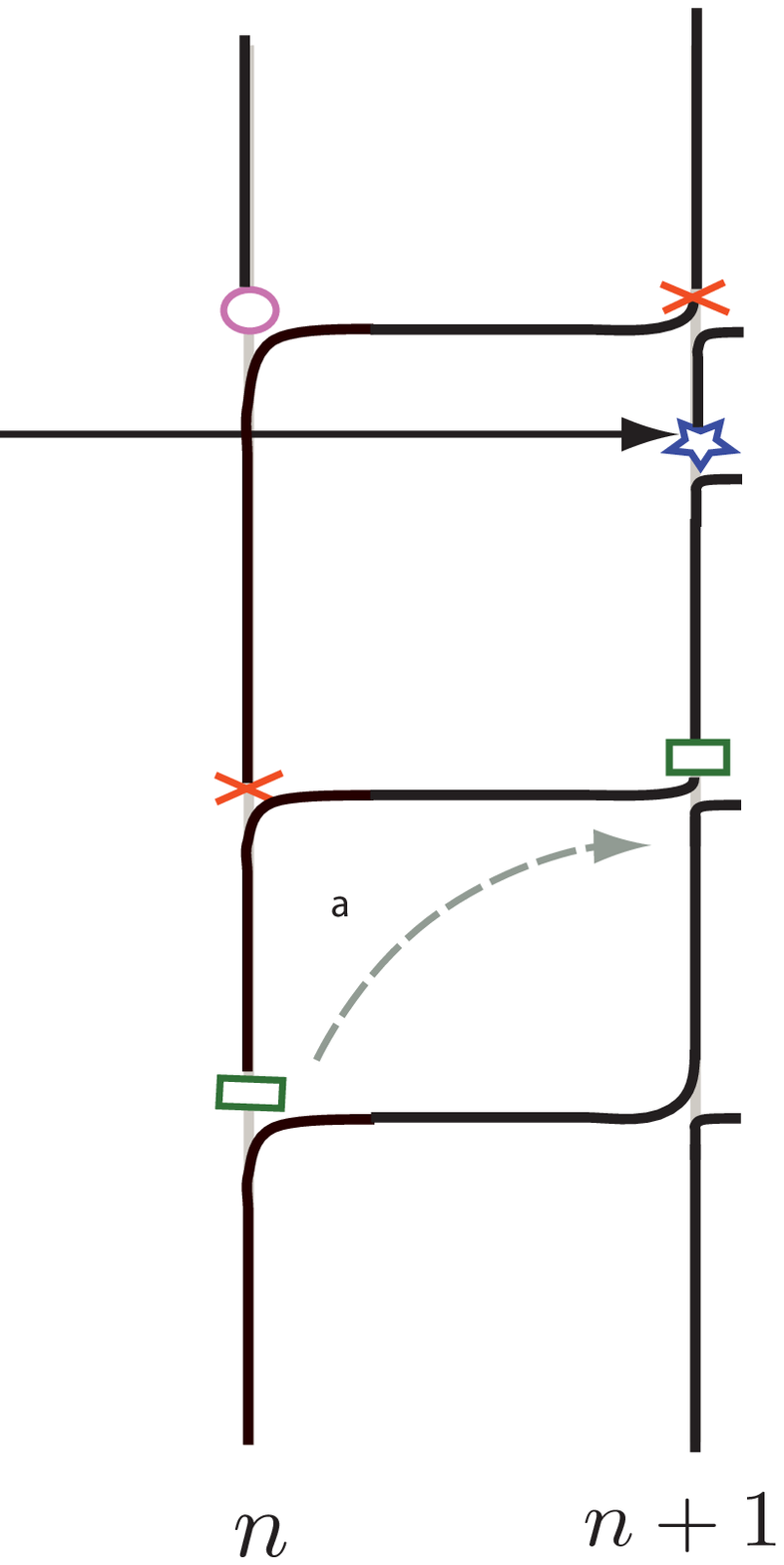}
\captionof {figure}{}
\label{Fig_eta}
\end{center}
\end{minipage}
\end{flushleft}
$\quad$ To each point of  $\eta_{2, \dots, n}^n$ we associate a label which records the level at which the line arriving at this point was born. By the induction assumption these labels are independent and take value $k$ with probability $(k-1)/ {\binom n 2}$. The birth level of a line arriving at time $t$ at level $n+1$ is the birth level of the line arriving at time $\phi^{-1}(t)$ at level $n$, and hence it is the label of the point $\phi^{-1}(t)$. The induction step is now completed by the following elementary observation: Consider an independent labeling of a stationary Poisson process $\eta=(\tau_i)_{i \in \mathbb{Z}}$, $\ldots <\tau_{-1}<\tau_0<\tau_1< \ldots$, on $\RR$ and perform an "upward shift" of this labeling, by assigning to each point $\tau_i$ as its new label the label of $\tau_{i-1}$. Then the new labeling has the same distribution as the old one.
\hfill$\square$

\newpage
{\it Proof of Theorem \ref{sec_th}}.

Let $s<t \in \RR$ be fixed. For every $k \ge 2$ we consider the sequence of lines born at level $k$ that die after time $s$, indexed by their death times $(t_i)_{i\ge1}$, with $s < t_1 < t_2 < \ldots$. For $i \in \mathbb N$ let $T_{ik}$ be the length of life of the $i$-th of these lines and let $M_k$ denote the number of these lines which die before time $t$. According to Proposition \ref{prop_poisson} the numbers $M_k$ are Poisson distributed with parameter $(k-1)(t-s)$ and independent from one another for different $k$'s. We show that
\[S_{s,t}:=  \sum_{k=2}^{\infty} \sum_{i=1}^{M_k} T_{ik}^2 \]
is infinite almost surely. 

To this end we first observe that for each $k \ge 2$ and $i \ge 1$ the random variable $T_{ik}$ has the same distribution as $T_k :=$  the sum of independent Exp${j \choose 2}$-distributed random variables, where $j$ ranges from $k$ to~$\infty$.

For each $k \geq 2$ we have
\begin{align}\notag
&\PP\Big(M_k  \notin \Big[\frac 1 2 (k-1)(t-s), 2(k-1)(t-s)\Big]
\\[0.5ex] \notag
&\quad\quad  \text{or} \quad\Big\{M_k \in \Big[\frac 1 2 (k-1)(t-s), 2(k-1)(t-s)\Big] \text{ and } T_{ik} \notin\Big[\frac 1 k, \frac 3 k\Big] \text{ for some } i=1,\dots, M_k \Big\} \Big) \notag
\\[0.5ex] \notag
& \leq \PP\Big(M_k  \notin \Big[\frac 1 2 (k-1)(t-s), 2(k-1)(t-s)\Big] \Big)
 \\[0.5ex] \notag 
 &\quad\quad  + \PP\Big( T_{ik} \notin\Big[\frac 1 k, \frac 3 k\Big] \text{ for some } i=1,\dots, \lceil 2 (k-1)(t-s)\rceil \Big\} \Big)
 \\[0.5ex] \label{BC}
& \leq \PP\Big(M_k  \notin \Big[\frac 1 2 (k-1)(t-s), 2(k-1)(t-s)\Big] \Big)+
 \lceil 2 (k-1)(t-s)\rceil \PP\Big(T_k \notin\Big[\frac 1 k, \frac 3 k\Big]\Big). 
\end{align}
Cram\'er's theorem guarantees that $\PP\Big(M_k  \notin \Big[\frac 1 2 (k-1)(t-s), 2(k-1)(t-s)\Big] \Big)$ decays exponentially in $k$ and hence the first term on the right-hand side is summable. For the second term we use Theorem 1 of \cite{DPS13} which says that the sequence $(kT_{k})_{k\ge2}$ (that converges a.s. to $2$ as $k\to \infty$) satisfies a large deviation principle with scale $k$ and a good rate function. Since
\[\PP\Big(T_{k} \notin\Big[\frac 1 k, \frac 3 k\Big]\Big) =\PP\Big(\Big|T_{k}-\frac 2 k\Big| > \frac 1 k\Big),\]
it follows that the second term on the right-hand side of (\ref{BC}) is also summable. By the Borel-Cantelli lemma we obtain that there exists an $\NN$-valued random variable $K_1 \geq 2$  such that for all $k\geq K_1$ 
\[M_k \in \Big[\frac 1 2 (k-1)(t-s), 2(k-1)(t-s)\Big] \quad \text{ and } \quad T_{ik} \in\Big[\frac 1 k, \frac 3 k\Big]  \quad \text{ for all } i=1, \dots,  M_k \quad  \text{ almost surely} \]
and in particular
\[M_k \geq \frac 1 2 (k-1)(t-s) \quad \text{ and } \quad T_{ik} \geq \frac 1 k  \quad \text{ for all } i=1, \dots,  M_k  \quad \text{ almost surely.}\]
Therefore, it holds that almost surely
\[
\sum_{k=K_1}^{\infty} \sum_{i=1}^{M_k} T_{ik}^2 \geq \sum_{k=K_1}^{\infty} \Big\lceil \frac 1 2 (k-1)(t-s)\Big\rceil \cdot \frac 1 {k^2}.
\]
Now since $K_1$ is almost surely finite, it follows that the sum on the right-hand side is infinite almost surely and that
\[S_{s,t} = \infty \quad \text{ almost surely},\]
which gives the claim.
\hfill$\square$

\section{Proof of Theorem~\ref{main_th}}
In order to prove Theorem~\ref{main_th} we first recall that Proposition~3.2 of \cite {PWW11} ensures the existence of a process $ \mathfrak L^{ld}$ having the same distribution as the Kingman tree length process  $\mathfrak L$ and such that $d_{Sk}(\mathfrak L^{ld,N}, \mathfrak L^{ld}) \to 0$ as $N \to \infty$ in probability. It thus suffices to prove  Theorem~\ref{main_th} for   $\mathfrak L^{ld}$ instead of $ \mathfrak L$.
\par
The following lemma is elementary; we include its proof for the sake of completeness.
\bl\label{old_lemma}
Let $(y_k)_{k \geq 1}$, $y_k: \RR \to \RR$  be a sequence of c\`adl\`ag functions satisfying that there exist two sequences $(\t_k)_{k \geq 1}$ and $(\g_k)_{k \geq 1}$ in $\RR$ such that $y_k$ has a jump of size $\g_k$ at time $\t_k$ for all $k \geq 1$. Moreover, suppose that the sequence $(y_k)_{k \geq 1}$ converges in the Skorohod topology to a c\`adl\`ag function $y$ and that the sequences $(\t_k)_{k \geq 1}$ and $(\g_k)_{k \geq 1}$ are convergent. Let $\t:=\lim_{k \to \infty} \t_k$ and $\g:=\lim_{k \to \infty} \g_k$ and assume that $\gamma \neq 0$. Then the function $y$ has a jump of size $\g$ at time $\t$.
\el
\bp
Let $\Lambda$ be the set of all strictly increasing and continuous  functions $ \lambda:[0,\infty] \to [0,\infty]$. Together with the stated assumptions, the convergence $d_{Sk}(y_k,y) \to 0$ implies the existence of a sequence $( \lambda_k )_{k \geq 1}$  of functions in $\Lambda$ such that 
\[ \rho_k:=\lambda_k(\tau_k) \to \tau \, \text{  and  } \,  \Delta y(\rho_k):=y(\rho_k)-y(\rho_k-)  \to \gamma \quad \text{ as } k \to \infty . \]
 If $\rho_k$ were different from $\tau$ for infinitely many $k$, then this would contradict the fact that large jumps of a c\`adl\`ag function are isolated.  Consequently, $\rho_k=\tau$ for all but finitely many $k$ (see also  \cite{JS} Proposition  VI.2.1 b) with $\alpha_k = \alpha = y$, $t_k=t=\tau$, $t_k' = \rho_k$).  Hence, $\Delta y(\tau) = \gamma$. \ep

\bpr\label{prop_2}
The sum of the squared jump sizes of the process $\mathfrak L^{ld}$ occurring in any interval of positive length is infinite almost surely.
\epr

\bp

Consider the look-down graph and recall that for every $N \in \NN$ the $N$-look-down graph can be recovered as the projection of the infinite graph onto its first $N$ levels. 

Let $\MG$ denote the set of all the lines in the infinite look-down graph and for a line $G=\bigcup_{j \in \NN_0}[s^G_j,s^G_{j+1}) \times \{l^G+j\}$ born at level $l^G$ let us set
\[
d^{G,N}:=\left\{ \begin {array} {ll} 
                  s^G_{N-l^G+1} \text{ , if } N \geq l^G \\ 
                  -\infty \text{ , otherwise,}
          \end {array} 
  			\right.
\]
the {\it exit time from level $N$} of the line $G$.

We are interested in the times and the sizes of the jumps of the processes $\mathfrak L^{ld,N}$. Jumps occur at the times $\{d^{G,N}\}_{G \in \MG}$ when lines die in the $N$-look-down process (i.e. they exit level $N$). Since for a fixed $G \in \MG$ the sequences $\{d^{G,N}\}_{N \in \NN}$ and $\{s^G_{N-l^G+1}\}_{N \in \NN}$ are identical for $N$ large enough, it follows that
\be\label{conv_d}
\lim_{N \to \infty} d^{G,N}=d^{G},
\ee
where $d^G$ is the death time of line $G$. The jump size of the process $\mathfrak L^{ld,N}$ at time $d^{G,N}$ has size equal to the life-length $T^{G,N}$ of the line $G$ up to level $N$ defined in (\ref{def_T_GN}). Note that the exponential times $X_j^G$ do not depend on $N$. Therefore, we have that
\be\label{conv_T}
\lim_{N \to \infty} T^{G,N}=\lim_{N \to \infty} \sum_{j=l^G}^{N}X_j^G=\sum_{j=l^G}^{\infty}X_j^G=T^G
\ee
almost surely, where $T^G$ defined in (\ref{def_T_G}) is the life length of line $G$ in the infinite look-down graph.
\par
In the following we fix an increasing sequence $(N_k)_{k \in \NN}$ in $\NN$ such that $d_{Sk}(\mathfrak L^{ld,N_k}, \mathfrak L^{ld})\to 0$ almost surely as $k \to \infty$. In view of (\ref{conv_d}) and (\ref{conv_T}) we now apply for every $G \in \MG$ Lemma~\ref{old_lemma} to the paths of $(\mathfrak L^{ld,N_k})_{k \geq 1}$, the sequence of times $(d^{G,N_k})_{k \geq 1}$ and the sequence of jump sizes $(T^{G,N_k})_{k \geq 1}$. Consequently,  for each $G \in \MG$, the limiting process $\mathfrak L^{ld}$ has a jump of size $T^G$ at time $d^G$. 

Thus, for the sum of the squared jump sizes of  $\mathfrak L^{ld}$ occurring in an interval $[0,t]$, $t>0$,
\[
\sum_{0 \leq s \leq t} (\Delta \mathfrak L^{ld}_s)^2 \geq \sum_{ G \in \MG: \atop d^G \in [0,t]  } \left(T^G\right)^2
\]
holds and since, according to Theorem~\ref{sec_th}, the right-hand side is infinite almost surely, the Proposition is proved.
\ep

\par
{\it Proof of Theorem~\ref{main_th}}

\par
It remains to show that for every $t>0$ any c\`adl\`ag path $X$  which obeys
$\sum_{0 \leq s \leq t} (\Delta X_s)^2=\infty$ 
has the property 
\be\label{eq_k}
\lim_{n \to\infty}\sum_{j=1}^{l^{(n)}} \Big(X_{\rho_j^{(n)}}-X_{\rho_{j-1}^{(n)}}\Big)^2 =\infty
\ee
for each sequence $(\mathcal P_n)_{n \in \NN}=\Big((\rho_j^{(n)})_{j=0,\dots,l^{(n)}}\Big)_{n \in \NN}$ of partitions of $[0,t]$ with mesh size tending to zero as $n \to \infty$. For this purpose we order the jump sizes of $X$ that occur in $(0,t)$ according to their sizes and denote by $(t_i)_{i\geq 1}$ the corresponding jump times, i.e. $|\Delta X_{t_1}| \geq |\Delta X_{t_2}| \geq \dots$ holds. Then, for every (fixed but arbitrary) $k \in \RR$ there exists an $m(k)$ such that
\be \nonumber
\sum_{i=1}^{m(k)} (\Delta X_{t_i})^2 \geq k.
\ee
For any jump time $t_i$ and every partition $\mathcal P_n$ let $\sigma_{i,n}$ be the largest point in the partition smaller than $t_i$ and $\tau_{i,n}$ be the smallest point in the partition larger than or equal to $t_i$. Then, for $n$ large enough, there is at most one of the $t_1,\ldots, t_{m(k)}$ between any two points of the partition $\mathcal P_n$ and thus
\[\sum_{j=1}^{l^{(n)}} \Big(X_{\rho_j^{(n)}}-X_{\rho_{j-1}^{(n)}}\Big)^2 \geq \sum_{i=1}^{m(k)} \Big(X_{\tau_{i,n}}-X_{\sigma_{i,n}}\Big)^2 \]
holds for $n$ large enough. 
Using the c\`adl\`ag property of $X$  we obtain that
\[\lim_{n \to \infty} \sum_{i=1}^{m(k)} \Big(X_{\tau_{i,n}}-X_{\sigma_{i,n}}\Big)^2 \ge \sum_{i=1}^{m(k)} {(\Delta X_{t_i})^2}.\]
Since $k$ was arbitrary, \eqref{eq_k} follows from the last three inequalities.

From this together with Proposition \ref{prop_2} and Lemma \ref{old_lemma}, Theorem~\ref {main_th} is immediate.
\hfill$\square$

\vspace{1cm}
\noindent
{\bf Acknowledgement.} We are grateful to Stephan Gufler, G\"otz Kersting and Etienne Pardoux for stimulating discussions. We also thank two referees for careful reading and helpful comments.

\end{document}